\documentclass[12pt]{article}
 \usepackage{amsrefs}
\usepackage{hyperref}
\usepackage{relsize}
\usepackage{enumitem}
\usepackage{graphicx}
\usepackage{float}
\usepackage{wrapfig}
\usepackage{picinpar}
\usepackage{tcolorbox}
\usepackage{hyperref}
\usepackage{epsfig}
 \usepackage{lmodern}
\usepackage{amsmath}
\def\bsq{\blacksquare}
\usepackage{tikz}

 \newif\ifstartedinmathmode
\newcommand\encircled[1]{%
  \relax\ifmmode\startedinmathmodetrue\else\startedinmathmodefalse\fi%
  \tikz[baseline,anchor=base]{%
  \node[draw,circle,outer sep=0pt,inner sep=.2ex]
    {\ifstartedinmathmode$#1$\else#1\fi};}%
}
\def\n{\noindent}

\usepackage{amsmath,amssymb,amsbsy,amsfonts,amsthm,latexsym,
            amsopn,amstext,amsxtra,euscript,amscd,amsthm}
\newtheorem{lem}{Lemma}
\newtheorem{lemma}[lem]{Lemma}

\newtheorem{thm}{Theorem}
\newtheorem{theorem}[thm]{Theorem}
\newtheorem{cor}{Corollary}
\newtheorem{corollary}[cor]{Corollary}

\def\\{\cr}
\def\({\left(}
\def\){\right)}
\def\[{\left[}
\def\]{\right]}
\def\<{\langle}
\def\>{\rangle}

\def\bsq{\blacksquare}
\def\eproof{$\hfill\bsq$\par}
\def\eq{\Leftrightarrow}
\usetikzlibrary{shapes}

\begin{document}

\title{Drainage Time and Shape: Inequalities from Torricelli’s Law}

\date{\today}

\pagenumbering{arabic}

\author{Eugen~ J. Ionascu    \footnote{Department of Mathematics, Columbus State University, Columbus, GA, USA Email: ionascu@columbusstate.edu }
}

\maketitle

\begin{abstract}We derive integral inequalities governing drainage time in convex solids, inspired by Torricelli’s Law, and introduce the Torricelli number as a shape invariant.
We use these considerations to construct a class of solids that can be used in building asymmetrical clepsydrae.
\end{abstract}



\maketitle

\section{Introduction}
We are interested in various consequences of the Torricelli’s Law of drainage of a liquid contained in a solid through a small hole situated at the bottom of the solid.  Evangelista Torricelli (1608-1648), who was Galileo’s last student \cite{Groetsch08},  formulated this law by saying that the speed flowing out of an orifice under gravity is equal to the speed a body would acquire falling freely from the same height: $v = \sqrt{2gh}$.  Torricelli didn't have calculus concepts, but these days, this is a standard topic usually in a differential equations course in the undergraduate curriculum for engineers/math/science majors since the law can be written as a simple autonomous separable differential equation

\begin{equation}\label{eq1} A(h)\frac{dh}{dt}=-k a\sqrt{2gh}
\end{equation} 
\noindent where $a$ is the area of the hole,  $A(h)$ is the area of the section at height $h$, $\frac{dh}{dt}$ is the instantaneous change in height of the level of liquid (vertical velocity considered with positive orientation upwards- hence the negative sign on the right-hand side of the equation), $k$ is a constant which depends of the liquid viscosity and some other factors such as friction along the exit hole, $g$ is the constant of acceleration due to gravity, and $h$ is the height of the water at time $t$ (arbitrary point in time). We will refer to (\ref{eq1}) as the {\it  Torricelli’s differential equation law}.

$$\underset{Figure\ 1}{\epsfig{file=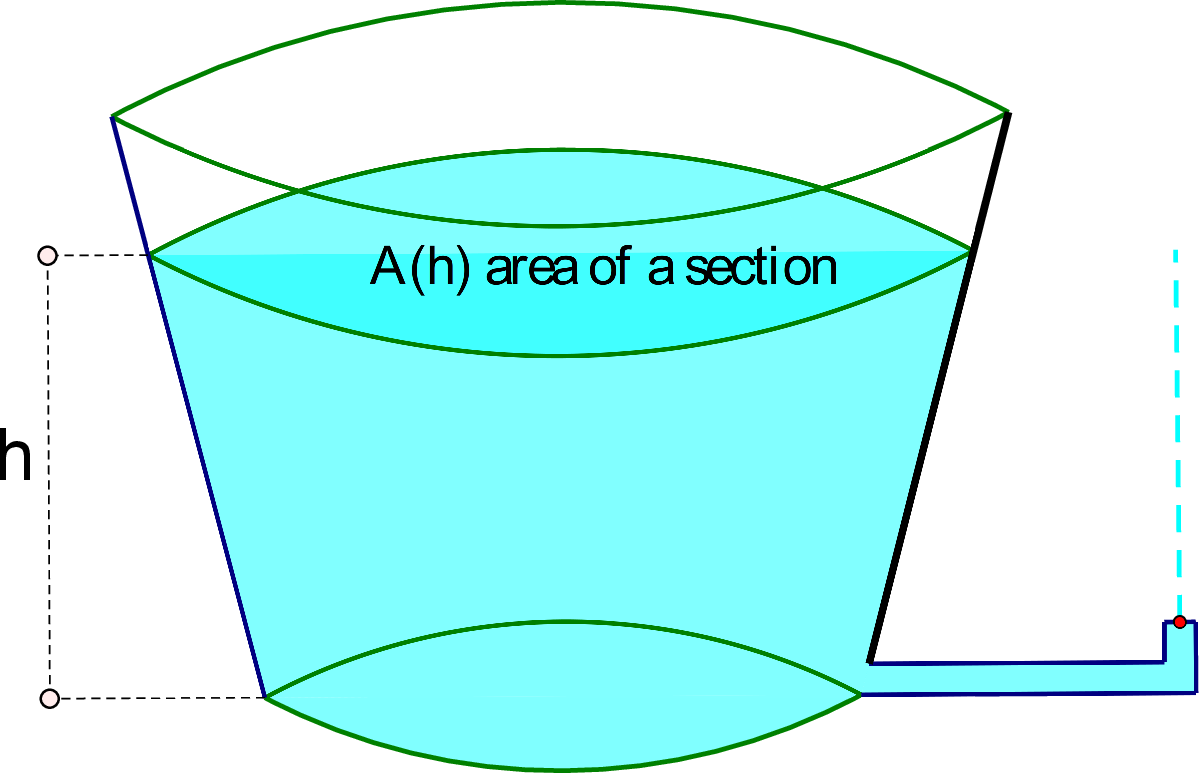,height=1in,width=2in}}$$

\n The key idea in the derivation of the law, was the experiment done (see Figure~1), in which it was discovered that the liquid stream has enough velocity at the outlet to push it to quite the same level as the level of the liquid in the recipient (under ideal setings, disregarding air resistance, and losses of energy through friction, etc).

We apply (\ref{eq1})  to certain types of solids and compare various times related to the same solid but maintaining the constants involved ($a$, $k$, and  $g$) in (\ref{eq1}).  We raise some questions about these times and solve some of them. 

\vspace{0.1in} 

To be a little more specific, we will assume that $h(0)=H$ and $h(T)=0$ for some positive value in time $T$-{\it  the total drainage time of the liquid for a certain solid and orientation}. Assuming a precise formula of $A$ in terms of $h$ is known, we can integrate (\ref{eq1}) over the interval $[0,T]$. In other words, we get 

$$\int_0^T  A(h)h^{-1/2} \frac{dh}{dt}dt= -ka\sqrt{2g}\int_0^Tdt = -ka\sqrt{2g}T.$$

\n The function $\frac{A(h)}{\sqrt{h}}$ is undefined and unbounded (in most cases) at/around  $h=0$ but the Riemann integral is convergent at this point since $$\int_t ^1 \frac{1}{\sqrt{h}}dh=2\sqrt{h}|_t^1=2(1-\sqrt{t})\ \ \ \to 2\ \ \ \text{as}\ \  t\to 0 \ (t>0).$$
\n So, we can define $F(x)=\int_0^x \frac{A(h)}{\sqrt{h}}dh $ and obtain a continuous function, given a continuous function $A$, with $F(0)=0$. With this, the equation above reduces to (by a change of variables $x=h(t)$, $t\in [0,T]$)
  
$$F(0)-F(H)= -ka\sqrt{2g}T\  \implies \ T=\frac{F(H)}{ka\sqrt{2g}}.$$

\n For a certain solid we will assume that $K=ka\sqrt{2g}$ is the same for all situations. In practice, this amounts to having the material  
 that bounds the solid being homogenous in thickness and $a$ (the hole area) being the same. As a result, we will record the formula above by 
 
$$T=\frac{F(H)}{K}=\frac{\int_0^H \frac{A(h)}{\sqrt{h}}dh}{K}=\frac{2 \int_0^{\sqrt{H}}A(s^2)ds}{K} \ \ \text{or}\ \ $$

\begin{equation}\label{eq2}
 \boxed{T=\frac{1}{K}\underset{\cal S}{\iiint}  \frac{1}{\sqrt{h}}dxdy dh  =\frac{1}{K}\underset{\cal S}{\iiint}  \frac{1}{\sqrt{h}}dhdxdy  =\frac{2}{K}  \int_0^{\sqrt{H}}A(s^2)ds, }
\end{equation} 
\n where $\cal S$ is our three-dimensional solid. The last expression makes it clear that $T$ is well defined for any non-negative continuous function $A$ on some interval $[0,H]$.

\vspace{0.1in}  A few comments about the last expression in (\ref{eq2}) are in order. Let us assume that our solid is like a coconut (not a perfect sphere in shape but convex) and we want to place the hole in a spot that makes $T$ the smallest possible (the orientation in space is given by having the coconut ``tangent" to the horizontal plane at the hole). How do we go about it?  What about a similar question but about a simpler solid like a platonic solid? Where is the best place to make an identical hole, assuming uniform material thickness (very small), for a regular tetrahedron for example?  

What if we are looking to maximize $T$?  Given a certain solid, we can define 
\begin{equation}\label{eqMAIN}\rho_{torr}= \frac{T_{max}}{T_{min}}=\frac{\underset{\cal S}{\iiint}  \frac{1}{\sqrt{h_{max}}}dv}{\underset{\cal S}{\iiint}  \frac{1}{\sqrt{h_{min}}}dv} =\frac{\int_0^{\sqrt{H_{max}}}A_{max} (s^2)ds}{\int_0^{\sqrt{H_{min} }}A_{min} (s^2)ds}
\end{equation} 
\n  which is independent of $K$, a mathematical constant that depends only on the shape of the solid, and as we can see it is independent of the drainage physical experiment. Let us call this number the {\bf Torricelli number} associated with a solid. For sure, a sphere as expected must have $\rho_{torr}=1$. What is $\rho_{torr}$ for an ellipsoid 
$${\cal E}=\{(x,y,z) | \frac{x^2}{u^2}+\frac{y^2}{v^2}+\frac{z^2}{w^2}=1\} , $$

\n a cylinder, or even for a box $[0,u]\times[0,v]\times [0,w]$?  How big can $\rho_{torr}$ be? Also, if $\rho_{torr} =1$ should we expect a spherical solid? In this paper, we are going address these questions for some simple solids. Some estimates of $\rho_{torr}$ are given in Section~2 for solids with a central symmetry. From a physics perspective, it is natural to think that the potential energy,
$$PE=\int_0^{H}A(h)h dh=\underset{\cal S}{\iiint} hdv$$ plays an important role in determining the maximum/minimum time given by (\ref{eq2}). More precisely, for bigger potential energy do we expect a shorter time?  In other words, if $h$ is maximized on average, isn't  $\frac{1}{\sqrt{h}}$ minimized on average? In Section~4, we will construct examples in which this is not the case. 

For many solids, the $T_{max}$ and $T_{min}$ are obtained simply by turning the solid upside down. For this reason, we compute such a ratio  (we call it {\it turn-up number}) in terms of a given orientation of the solid, which provides some information about  $\rho_{torr}$ above and in some instances it is precisely $\rho_{torr}$. We will discuss those in Section 3. In Section~4 we concentrate on solids of revolutions with the turn-up number equal to 1.

\vspace{0.2in} 

\section{Symmetrical solid with respect to a point.} 

\vspace{0.2in} We are interested in  convex solids $S$ that are having the property that there exits a point $C$ (the center) such that 
for every point $P$ of $S$, the reflection of $P$ with respect to $C$ say $P'=R(P)$ is also in $S$, i.e., the vectors $\overrightarrow{CP}$ and $\overrightarrow{CR(P)}$ are opposite to one another ($\overrightarrow{CP}+\overrightarrow{CR(P)}=0$).
Examples of such solids abound. For instance,  a sphere, a box, an ellipsoid, all the platonic solids, any solid of revolution determined by a symmetric function, etc.  

$$\underset{Figure\ 2}{\epsfig{file=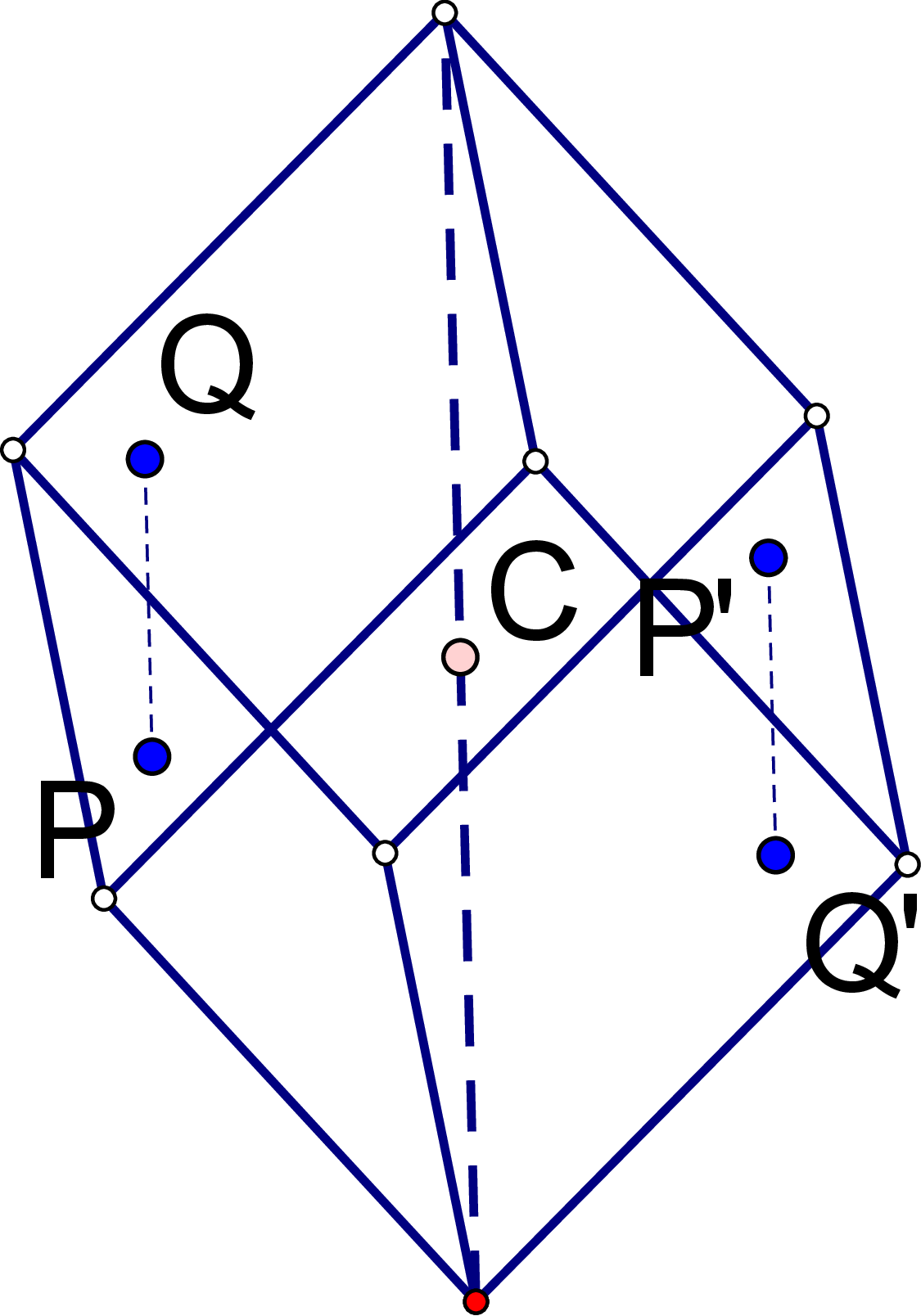,height=2in,width=1.5in}}$$

Let us assume that $S$ (we will think of a cube, as in Figure~2) is in an arbitrary position and the orifice is at one of the vertices of the cube which insures the drainage of all the liquid. Then by formula (\ref{eq2}) we have 

\begin{equation}\label{eq2ineq}
T=\frac{\int_0^H \frac{A(h)}{\sqrt{h}}dh  }{K} =\frac{1}{K}\int_0^H \underset{Region(h)}{\int \int} \frac{1}{\sqrt{h}} dxdy dh
\end{equation} 
\n where $Region(h)$ is the intersection of a plane at hight $h$ with the cube which has area $A(h)$. We can use Fubini's Theorem and interchange the order of integration (we are dealing with a non-negative integrant which is continuous of all three variables). We can then rewrite the above time as 

$$T=\frac{1}{K}\underset{\cal R}{\int\int} \int_{P(x,y)}^{Q(x,y)}\frac{1}{\sqrt{h}}dh  dxdy, $$

\n where $P(x,y)$ and $Q(x,y)$ are the points where the vertical on top of $(x,y)$ intersects the surface of the cube, and $\cal R$ is the region in the xy-plane 
for which such interesctions exist. 
Because of the symmetry, we have 
$$2T=\frac{1}{K}\underset{R}{\int\int}( \int_{P(x,y)}^{Q(x,y)}\frac{1}{\sqrt{h}}dh  dxdy+\int_{Q'(x,y)}^{P'(x,y)}\frac{1}{\sqrt{h}}dh  dxdy), $$
\n with the notation discussed above and shown in Figure~2.  Now the integration with respect to $h$ is $2\sqrt{h}$ and we can continue
$$T=\frac{1}{K}\underset{R}{\int\int}(\sqrt{h_Q}-\sqrt{h_P}+\sqrt{h_{P'}}-\sqrt{h_{Q'}})dxdy.$$
Now, we observe that $\sqrt{a}-\sqrt{b}=\frac{a-b}{\sqrt{a}+\sqrt{b}}$ and $\frac{\sqrt{a}+\sqrt{b}}{2}\le \sqrt{\frac{a+b}{2}}$. This implies that
 $$\sqrt{a}-\sqrt{b}=\frac{a-b}{\sqrt{a}+\sqrt{b}}\ge \frac{a-b}{2\sqrt{\frac{a+b}{2}}   }.$$
We apply this inequality for $h_Q$ and $h_{Q'}$ and obtain that 
$$\sqrt{h_Q}-\sqrt{h_{Q'}} \ge \frac{h_Q-h_{Q'}}{2\sqrt{h_C}}$$ since $\frac{h_Q+h_{Q'}}{2}=h_C$.
We have a similar inequality for $\sqrt{h_{P'}}-\sqrt{h_{P}}$. Adding them together and integrating, the formula above for $T$ turns into an inequality:
$$T\ge \frac{1}{2K\sqrt{h_C}}\underset{R}{\int\int}(h_Q-h_P+h_{P'}-h_{Q'})dxdy= \frac{1}{2K\sqrt{h_C}}\cdot 2V= \frac{V}{K\sqrt{h_C}},$$
\n where $V$ is the volume of the solid $S$. So, we have shown that $T\ge  \frac{V}{K\sqrt{h_C}}$.

\par \vspace{0.1in} On the other hand, we observe that if $a+b=c$ then we also have the inequality $\sqrt{a}-\sqrt{b}\ge \frac{a-b}{\sqrt{c}}$, $a>b\ge 0$. This is equivalent to $\sqrt{c} \ge \sqrt{a}+\sqrt{b}$ which is clearly true if we square both sides. Using the same arguments as before but using this inequality we obtain $T\le \frac{2V}{K\sqrt{2h_C}}$. Therefore, we have the following estimates of the time of drainage in the case of a symmetric solid.

\begin{theorem}\label{symmineq}
For a symmetric, convex solid $S$ with respect to a point $C$ the time of drainage $T$ through an orifice at the bottom of the solid   
located in the $xy$-plane satisfies the inequalities
\begin{equation}\label{symmetric} 
\frac{V}{K\sqrt{h_C}}\le T\le \frac{\sqrt{2}V}{K\sqrt{h_C}},
\end{equation}
\n where $h_c$ is the altitude of $C$ relative to the $xy$-plane.
\end{theorem}

We then obtain an estimate for the Torricelli's number. For a symmetric solid let us define $D=max\{|PP'||P\in S\}$ and $d=max\{|PP'||P\in S\}$.

\begin{corollary}\label{torrineq}
For a symmetric, convex solid $S$ with respect to a point $C$ the  Torricelli's number satisfies 
\begin{equation}\label{torrsymmetric} 
\rho_{torr}\le \sqrt{2}\sqrt{\frac{h_{c\ max}}{h_{c\ min}}}=\sqrt{\frac{2D}{d}},
\end{equation}
\n where $h_{c\ max}$  and $h_{c\ min}$ are the maximum/minimum possible altitude of $C$ relative to the orifice plane.
\end{corollary}

The following fact is expected but we have a simple proof for it.

\begin{corollary}\label{torrunbounded}
 The Torricelli number is unbounded.
\end{corollary}

\n {\bf Proof:} Let's consider a box of the size $S=1\times 1\times n$ with $n$ big. First, we position $S$ on one of the faces so that the long dimension is along the horizontal. The center is then at a distance of $1/2$ to the $xy$-plain. We obtain the inequalities from (\ref{symmetric}) 
$$  \frac{n\sqrt{2}}{K}\le T_1\le \frac{2n}{K}.$$ 
Then we position $S$ with the long dimension vertically. The center is then at a distance of $n/2$ to the $xy$-plain. Again, from (\ref{symmetric}) 
$$\frac{\sqrt{2n}}{K}\le T_2\le \frac{2\sqrt{n}V}{K}.$$
Then, $T_{min}\le T_2\le \frac{2\sqrt{n}}{K}$ and $T_{max}\ge  T_1\ge \frac{n\sqrt{2}}{K}$. Then, we get that 
$$\rho_{torr}\ge \frac{n\sqrt{2}}{2\sqrt{n}}=\frac{\sqrt{n}}{\sqrt{2}}\to \infty,$$
\n which implies that $\rho_{torr}$ is unbounded. \eproof

\subsection{Cubes}

 For a cube, we get $\rho_{torr}\le \sqrt{2\sqrt{3}}$. It is clear from (\ref{symmetric})  that $T_{max}= \frac{\sqrt{2}V}{K\sqrt{h_C}}= \frac{2V}{K}$ which is attained when the cube is positioned with two of its opposite faces horizontally. For  $T_{min}$ we will show that it is attained in the position shown in Figure~2, with one of its big diagonals being vertical. It is a non-trivial exercise to show that in this case (for a cube of sidelengths 1), the area function is 
$$A(h)=\begin{cases} \frac{3(\sqrt{3}-h)^2\sqrt{3}}{2} \ \ \text{if}\ \ h\in [\frac{2\sqrt{3}}{3},\sqrt{3}], \\ \\ 
\frac{3\sqrt{3}(2\sqrt{3}h - 1 - 2h^2)}{2}\ \ \text{if}\ \ h\in [\frac{\sqrt{3}}{3},\frac{2\sqrt{3}}{3}], \\ \\
\frac{3h^2\sqrt{3}}{2} \ \ \text{if}\ \ h\in [0,\frac{\sqrt{3}}{3}].
\end{cases} 
$$
One can check that $\int_0^{\sqrt{3}}A(h)dh=1$ (the volume) and 

\begin{equation}\label{mintforcube}
T=\frac{1}{K}\int_0^{\sqrt{3}}\frac{A(h)}{\sqrt{h}}dh=\frac{8\sqrt[4]{3}}{5K} \left(1+3\sqrt{3}-4\sqrt{2}\right).
\end{equation}

We will show that this is $T_{min}$ and so for a cube we have

\begin{equation}\label{torrforcube}
\rho_{torr}=\frac{5}{4\sqrt[4]{3}\left(1+3\sqrt{3}-4\sqrt{2}\right)}\approx 1.7611678552583516780\ .
\end{equation}

In particular, we see that the inequality $\rho_{torr}\le \sqrt{2\sqrt{3}}$ is a good estimate since $\sqrt{2\sqrt{3}}\approx 1.861209718$.
Another interesting fact here is that $\int_0^{\sqrt{3}}A(h)hdh=\frac{\sqrt{3}}{2}$, which means the expected result that the potential energy is exactly that of a mass equal to the mass of the liquid located at the center of mass. The graph of $A$ is included in Figure~3. 
$$\underset{Figure\ 3}{\epsfig{file=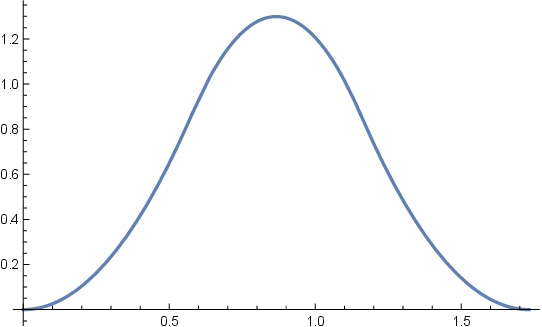,height=1in,width=2in}}$$
We can treat $A$ as a probability distribution. So, for example, we can compute the probability that a random point in the cube is located at an altitude 
$h$ between $h_1$ and $h_2$:   $\int_{h_1}^{h_2}A(h)dh$. Despite our complicated expression for $A$, it looks like this probability is a rational number for every $h_1/\sqrt{3}$ and $h_2/\sqrt{3}$ rational numbers. As an example, we get that 
$$ \int_{0}^{5\sqrt{3}/9}A(h)dh=\frac{101}{162}.$$

We will need the following optimization result.
\begin{lemma} Let us assume that we have two distinct points on the unit circle: $P(\cos(t),\sin(t))$ and $Q(\cos(s),\sin(s))$ in the positive half-plane ($y>0$). We consider the function
$$f(u)=\frac{1}{\sqrt{\sin(t+u)}}+\frac{1}{\sqrt{\sin(s+u)}}$$
\n defined for all $u$ so that $\sin(t+u)>0$ and $\sin(s+u)>0$. Then $f$ has a minimum which is attained at $u_{min}=\frac{\pi}{2}-\frac{s+t}{2}$, i.e., $P$ and $Q$ rotated with angle $u$ gives a segment $P'Q'$ parallel to the $x$-axis. This minimum is the only critical point in the interval $[0,2\pi]$, so a movement of $PQ$ decreases if $u$ is increasing but smaller than $u_{min}$. 
\end{lemma}

\n {\bf Proof:} Let us assume that $s$ and $t$ are in the $[0,\pi]$. Taking the derivative of $f$ we obtain 
$$f'(u)=-\frac{\cos (s+u)}{2 \sin ^{\frac{3}{2}}(s+u)}-\frac{\cos (t+u)}{2 \sin ^{\frac{3}{2}}(t+u)}.$$

\n To look for critical points we need first to require that 
$$\cos (s+u)\cos (t+u)=\frac{1}{2}[\cos (s+t+2u)+\cos (s-t)]<0.$$
We observe that this condition is satisfied if $s+t+2u=\pi$ because $\cos \pi=-1$ and $\cos (s-t)$ is not equal to $1$. So, the equation for critical points is then equivalent, under this constrain, with 
$$\cos^2 (s+u) \sin ^3(t+u)=\cos^2 (t+u) \sin ^3(s+u).$$
This is equivalent to 
$$[1-\sin^2 (s+u)] \sin ^3(t+u)=[1-\sin^2 (t+u)]\sin ^3(s+u),$$

\n and if we set $\sin(t+u)=x$ and $\sin(s+u)=y$, we get an equation only in $x$ and $y$:

$$x^3-x^3y^2=y^3-y^3x^2 \ \  \eq\ \  x^3-y^3+x^2y^2(y-x)=0\ \ \eq \ \ (x-y)(x^2+xy+y^2-x^2y^2)=0.$$

Clearly $x^2+xy+y^2-x^2y^2=x^2+xy+(1-x^2)y^2>0$ since $x,y>0$ and $1-x^2\ge 0$. Then, there is only one option, which is $x=y$. In this case, either $t+u=s+u+2k\pi$ ($k$ integer) which is not possible for  $t,s\in [0,\pi]$. The other alternative is $t+u=-(s+u)+(2k+1)\pi$ ($k$ integer) and of course, we are interested in the smallest solution in terms of $u$ which is exactly what we predicted  $u=\frac{\pi}{2}-\frac{s+t}{2}$. Adding a $\pi$ to $u$ makes $f$ undefined so we have to add a multiple of $2\pi$ which is not of interest since $f$ is  $2\pi$-periodic. So, this critical point is the only critical point and because the function is unbounded from above this point can only be a point of minimum.  \eproof

This optimization result can be generalized in the following way.

\begin{lemma} Let us assume that we have two distinct points on the unit circle: $P(\cos(t),\sin(t))$ and $Q(\cos(s),\sin(s))$. We consider the function
$$f(u)=\frac{1}{\sqrt{a\cos(t+u)+b\sin(t+u)+c}}+\frac{1}{\sqrt{a\cos(s+u)+b\sin(s+u)+c}}$$
\n defined for all $u$ so that the square roots make sense. Then $f$ has a minimum which is attained at $u_{min}$ that makes the two terms in $f(u)$ equal.     This minimum is the only critical point in the interval $[0,2\pi]$, so a movement of $PQ$ decreases if $u$ is increasing but smaller than $u_{min}$. 
\end{lemma}

We leave the proof for the reader. 

\begin{theorem}\label{mintimecube}
The time $T$ computed above in (\ref{mintforcube}) is $T_{min}$ for the cube. 
 \end{theorem}

$$\underset{Figure\ 4}{\epsfig{file=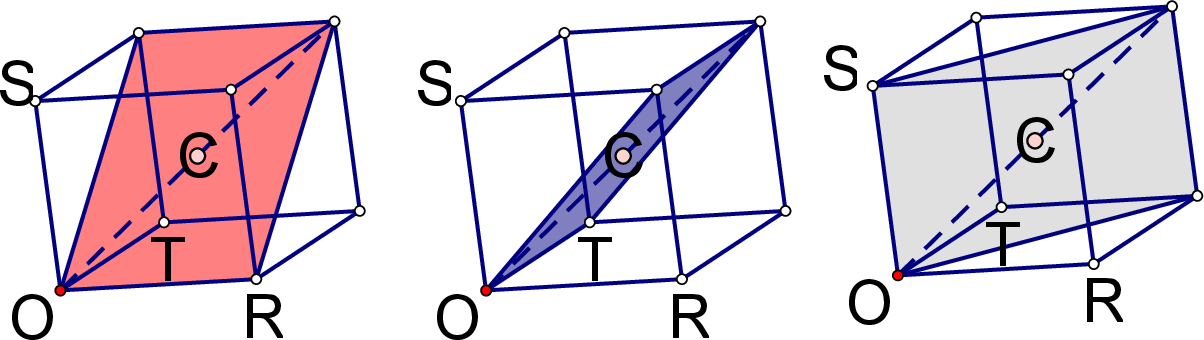,height=1.5in,width=3.3in}}$$ 
\n {\bf Proof:}  To show that $T$ computed above is $T_{min}$, let us start with a cube in an arbitrary position with the orifice on the bottom in one of the vertices of the cube say vertex $O$. We want to show that in two rotations around $\vec{OR}$, $\vec{OT}$, or $\vec{OS}$ (Figure~4) the cube can end in the position with a big diagonal vertical. For each such rotation the average value of $\frac{1}{\sqrt{h}}$ is going down which shows the integral over the cube is going down and so the result follows. 

The first thing to observe is that in Figure~4 we colored the three possible planes that divide the cube into two symmetrical parts, i.e., each part is the reflection of the other within that plane.  The rotations need to be done in such a way that these planes move toward a vertical position in order to decrease the average of $\frac{1}{\sqrt{h}}$. Let's write $\vec{k}$ (the unit vector that gives the vertical direction) in terms of  $\vec{OR}$, $\vec{OT}$ or $\vec{OS}$:
$$\vec{k}=r\vec{OR}+t\vec{OT}+s\vec{OS}.$$
We are working with the unit cube so we know that $r^2+t^2+s^2=1$. Since the cube is supposed to have all vertices above the xy-plane, we get that $r$, $t$ and $s$ are positive or zero. Not all these numbers are more than $\frac{1}{\sqrt{3}}$. Then at least one is less than or equal to  $\frac{1}{\sqrt{3}}$. Say, that number is $s$. If it is equal to  $\frac{1}{\sqrt{3}}$, we skip the first rotation. So let us assume that $r<\frac{1}{\sqrt{3}}$. A similar argument can be used to conclude that the other two numbers cannot be all less or equal than $\frac{1}{\sqrt{3}}$, so let us assume that $s>\frac{1}{\sqrt{3}}$. We will do a rotation around $\vec{OT}$, so that $r'$ becomes $\frac{1}{\sqrt{3}}$.

The big diagonal $OO'$ as a vector is $\vec{OO'}=\vec{OR}+\vec{OT}+\vec{OS}$, and the plane that rotates is spaned by $\vec{OT}$ and $\vec{OR}+\vec{OS}$. The projection of $\vec{k}$ on this plane is 
$$t \vec{OT}+\frac{r+s}{\sqrt{2}}(\vec{OR}+\vec{OS})$$
\n and its norm is $\sqrt{t^2+(r+s)^2/2}$. So, the plane moves toward the vertical position if $r+s$ gets bigger. From the identity $(r+s)^2+(s-r)^2=2(r^2+s^2)$ since we know $r^2+s^2=1-t^2$, along the rotation is constant. So, $r+s$ gets bigger if and only if $s-r$ gets smaller, or equivalently $r$ gets bigger. So, we will do a rotation until $r'=\frac{1}{\sqrt{3}}$. Once we do that we repeat the argument around $\vec{OR'}$ and then obtain that $\vec{OO'}$ becomes vertical. \eproof

\vspace{0.1in}
Another approach to calculating these times is to use Gauss' formula:

$${\displaystyle \iiint \limits _{V}\left(\mathbf {\nabla } \cdot \mathbf {F} \right)\,\mathrm {d} V=}  \iint \limits _{\partial V}  {\displaystyle (\mathbf {F} \cdot \mathbf {\hat {n}} )\,\mathrm {d} S}$$
\n by taking the vector field $F=<0,0,2\sqrt{z}>$. Then, the formula (\ref{eq2}) becomes

\begin{equation}\label{eq2bis}
 \boxed{T=\frac{1}{K}\underset{\cal S}{\iiint}  \frac{1}{\sqrt{z}}dxdydz   =\frac{2}{K} \iint \limits _{\partial S}\sqrt{z}\hat {n}_z \mathrm {d} S, }
\end{equation} 

\n where $\hat {n}_z$ is the third component of the exterior normal unit vector to the surface $\partial S$, i.e.,  $\hat {n}=<\hat {n}_x,\hat {n}_y,\hat {n}_z>$. The advantage of using such a formula is that this unit evector is constant for each face of the polyhedron $S$ if we are dealing with such a solid.

\vspace{0.2in}
\subsection{Regular Octahedrons}

For the case of the octahedron, let us assume without loss of generality that the side lengths are equal to $1$. 

$$\underset{Figure\ 5(a),\  T_{max} }{\epsfig{file=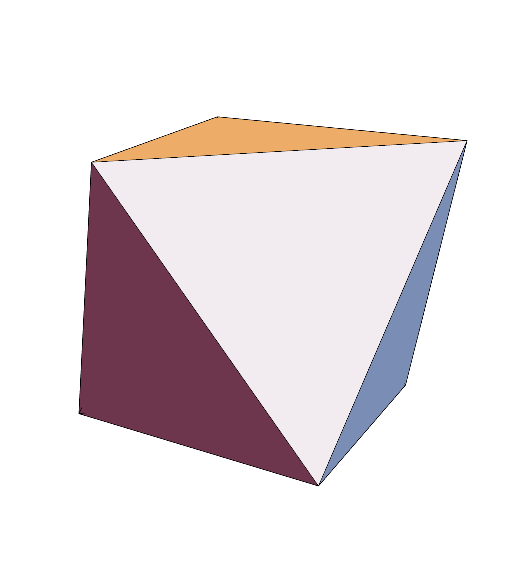,height=2in,width=2in}}\underset{Figure\ 5(b),\  T_{min} }{\epsfig{file=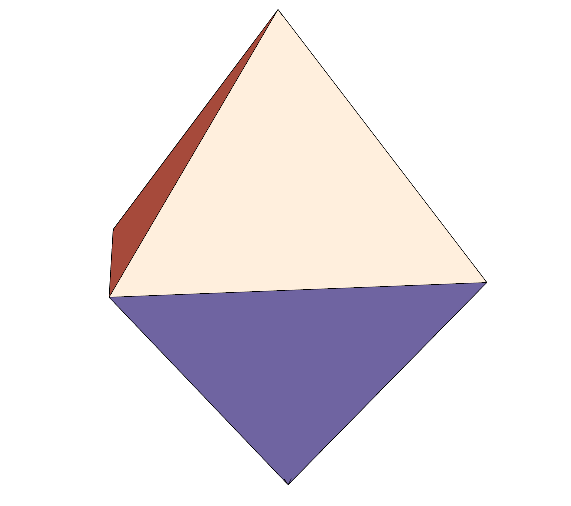,height=2in,width=2in}}$$
For each vertex of the octahedron, we have two symmetry planes containing an edge incident to that vertex. In fact, each one of these contains two of the edges incident to the specific vertex. We can use these two planes to show that the extreme positions are the ones in Figure~5, the same way as in the case of the cube. These planes become vertical for the minimum position. In the maximum position both of these planes make congruent angles with the horizontal plane, which are half of the dihedral angle of the octahedron: $\arccos(1/\sqrt{3})\approx 54.73561^{\circ}$. 
 
First, for the maximum time we get the following area function when calculating the maximum drainage time:
$$A_{max} (h)= \begin{cases}
\frac{a}{4}((2 - \frac{h}{H})^2 - 3(1 - \frac{h}{H})^2) \ \ \text{if} \ \ h \in [\frac{H}{2}, H], \\
\frac{a}{4}((1 + \frac{h}{H})^2 - 3( \frac{h}{H})^2) \ \ \text{if} \ \ h \in [0,\frac{H}{2}) \\
\end{cases}=\frac{a}{4}(1+\sqrt{6}h-3h^2), $$ 
where $a = \sqrt{3}, b = \sqrt{2},$ and $H=H_{max}=\frac{b}{a} =\frac{\sqrt{6}}{3}$. 

  Using this formula we calculate the maximum time
$$T_{max}=\frac{1}{K}\frac{19}{5\cdot 6^{3/4}}$$

As a checking, to get the volume we simply integrate $A(h)$, which gives us $V=\frac{\sqrt{2}}{3}$. For another checking, one can see that $$\frac{ \int_0^H A(h) \cdot hdh }{V}= \frac{H}{2}$$

\n as expected  (because of the symmetry).  As in the case of the cube, we also notice that cutting the octahedron with a horizontal plane through rational fraction of $H$ gives a rational fraction of $V$:

$$\frac{\int_0^{tH} A(h)dh}{V}=\frac{3t+3t^2-2 t^3}{4}\in \mathbb Q, \ \ t\in \mathbb Q.$$

Then, we figure out the area function when calculating the minimum drainage time (Figure 5(b)):
$$A_{min} (h)= \begin{cases}
2(H - h)^2 \ \ \text{if} \ \ h \in [\frac{H}{2}, H] \\
2h^2 \ \ \text{if} \ \ h \in [0,\frac{H}{2}) \\
\end{cases}$$ 
where $H =H_{min}= \sqrt{2}$.
  Using this formula for area,  we calculate the minimum time
$$T_{min}=\frac{1}{K} \cdot \frac{8}{15}(8 \cdot 2^{1/4}-5 \cdot 2^{3/4}).$$

Taking the ratio, and after some rationalizations we have
\begin{equation}\label{torrforpctahedron}
\rho_{torr}=\frac{19\sqrt[4]{3}(8+5\sqrt{2})}{224}\approx 1.68240255892043650...\ .
\end{equation}

\subsection{Regular Tetrahedrons}

$$\underset{Figure\ 6(a),\  T_{max} }{\epsfig{file=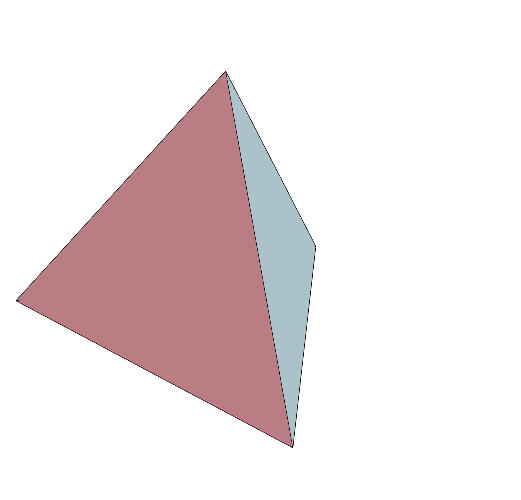,height=2in,width=2in}}\underset{Figure\ 6(b),\  T_{min} }{\epsfig{file=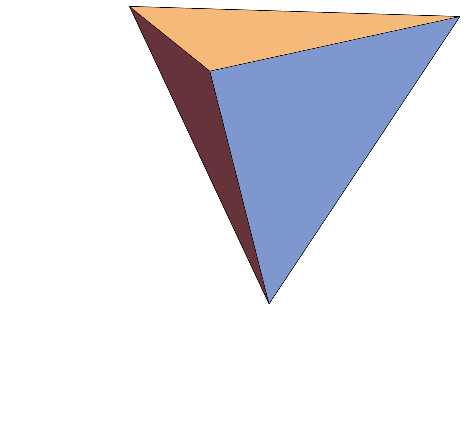,height=2in,width=2in}}$$
We will show that in this case 
\begin{equation}\label{torrfortetrahedron}
\rho_{torr}=\frac{8}{3}\approx 2.6666...\ .
\end{equation}
Let us work with side lengths of 1 again. Then the height is equal to $H=\sqrt{2/3}$. The formulae for the area in the maximum case is $A_{max}=\frac{3}{8}\sqrt{3}(H - h)^2$ and for the minimum is simply $A_{min}=\frac{3}{8}\sqrt{3}h^2$. Then we calculate $T_{max}=\frac{4 \sqrt[4]{2}}{5\cdot 3^{3/4}}$ and $T_{min}= \frac{\sqrt[4]{3}}{5\cdot  2^{3/4}}$ which give (\ref{torrfortetrahedron}). For each vertex, there are three symmetry planes that contain an incident edge to that vertex. The argument for extrema is similar to the previous ones. For the minimum these planes are vertical and for the maximum, two of the planes are making congruent  angles with the horizontal plane. These angles are the smallest possible since they are forced by the dihedral angles of the faces of the tetrahedron to
$\arctan(\frac{1}{\sqrt{2}})\approx 35.26439$.

\vspace{0.1in}

\subsection{Regular Icosahedrons}

 We will work with the icosahedron having vertices 
$$(0,\pm 1,\pm \phi), (\pm 1,\pm \phi, 0), (\pm \phi,0,\pm 1),$$
\n where $\phi=\frac{1+\sqrt{5}}{2}$, the Golden Ratio. This is the well-known construction of Luca Paccioli (see \cite{Joyce97}).  The edge lengths are then equal to $2$. It is known that the volume of this solid is 
${\cal V}=\frac{5}{6}\phi^2 s^3$  where $s$ is the length of each of its  edges, so in our case we get ${\cal V}=\frac{10}{3}(3+\sqrt{5})$. Let us introduce the notation $A=(\phi, 1, 0)$,  $B=(\phi, -1, 0)$, and $C=(1, 0, \phi)$. One can check that the equation of the plane determined by these three points is  $\phi x-(1-\phi) z=\phi^2$.

$$\underset{Figure\ 7(a),\  Paccioli's \ construction }{\epsfig{file=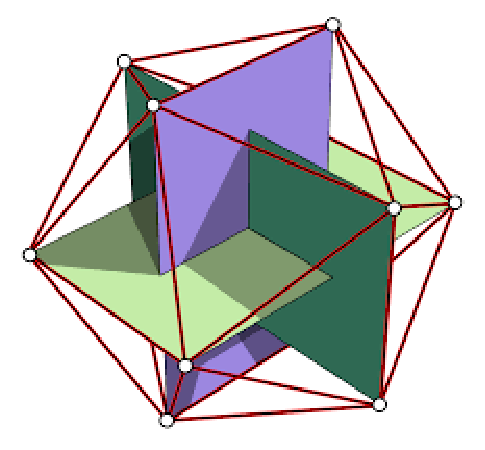,height=2in,width=2in}}
\underset{Figure\ 7(b),\  T_{min} }{\epsfig{file=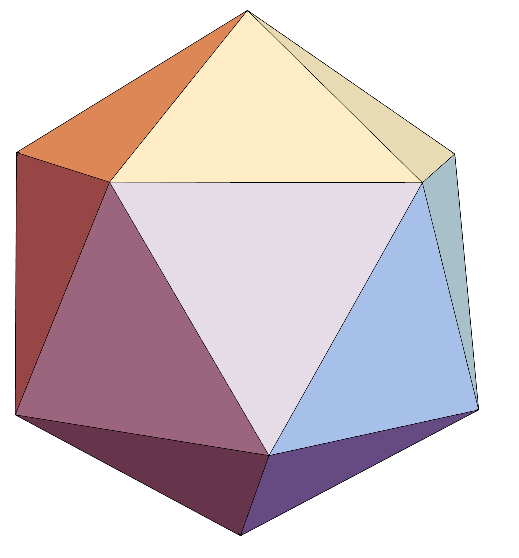,height=2in,width=2in}}
\underset{Figure\ 7(c),\  T_{max} }{\epsfig{file=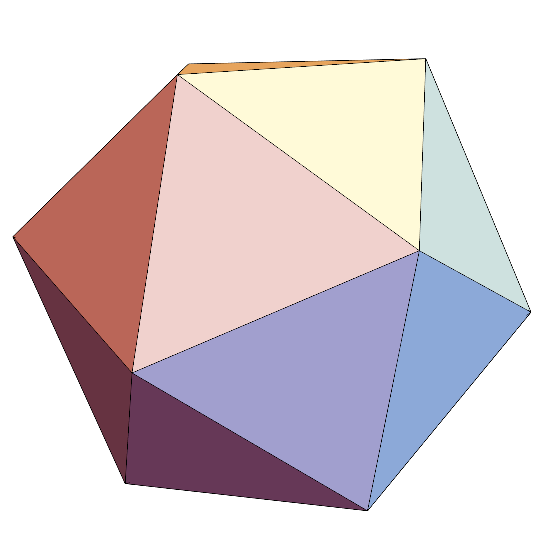,height=2in,width=2in}}$$

To compute the two times (see Figure~7), we need to solve Kepler's problem for this Platonic solid: ``find the radii of the inscribed and circumscribed spheres. " These radii give the two $H$-values. If $R$ is the radius of the circumscribed sphere, then $R=d(O,A)$ where $O$ is the origin of the coordinates system, $d$ is the Euclidean distance.  Then we get $$R=\sqrt{1+\phi^2}=\sqrt{1+\phi+1}=\sqrt{2+\phi}=\sqrt{\frac{5+\sqrt{5}}{2}}.$$
The radius $r$ of the sphere inscribed is the distance to the triangle $ABC$ from $O$ (the origin).    In general, the distance from a point
$P=(x_0,y_0,z_0)$ to a plane $ax+by+cz+d=0$ is equal to
$$\frac{|ax_0+by_0+cz_0+d|}{\sqrt{a^2+b^2+c^2}}.$$

 Then, we obtain
$r=\frac{\phi^2}{\sqrt{\phi^2+(1-\phi)^2}}$. Since $\phi^2=\phi+1$, we can simplify the expression for $r$ to
$$ r=\frac{\phi+1}{\sqrt{2\phi^2-2\phi+1}}=\frac{\phi+1}{\sqrt{3}}.$$
In the case of the minimum, $H_{min}=2R=2\sqrt{\frac{5+\sqrt{5}}{2}}$ (see Figure~7(b) ). For a regular pentagon of sides $s$, the radius of the circle circumscribed is 
$r'=\sqrt{\frac{5+\sqrt{5}}{10}}s$ and its area is ${\cal A}(s)=\frac{ \sqrt{5(5+2\sqrt{5})} }{4}s^2$. To determine the formula for $A_{min}$, we must find the splitting points in its piecewise definition: $h_s$ and $h'_s$ ($h_s<h'_s$). This leads to the height of the right pentagonal pyramid with base sides of $2$ and the oblique edges of 2. Pythagorean theorem gives $h_s=\sqrt{s^2-r'^2}=2\sqrt{\frac{5-\sqrt{5}}{10}}$. Then the ratio $\frac{h_s}{H_{min}}$ turns out to be equal to 
$$  \sqrt{\frac{5-\sqrt{5}}{5(5+\sqrt{5})}}=\frac{5-\sqrt{5}}{10}=\frac{1}{2}-\frac{\sqrt{5}}{10},$$

\n which gives us the other splitting point $h'_s$ which satisfies 
$$\frac{h'_s}{H_{min}}=\frac{1}{2}+\frac{\sqrt{5}}{10}.$$
The function $A_{min}(h)$, defined on the interval $[0,H_{min}]$,  is clearly symmetric with respect to the midpoint. For each $h\in [0,h_s]$, $A_{min}(h)$ is the area of a regular pentagon of size $s$, and $s$ depends linearly with $h$. In other words $A_{min}(h)=ch^2$ for some constant $c$, which can be determined from the equation 
$$A_{min}(h_s)=ch_s^2={\cal A}(2)=\sqrt{5(5+2\sqrt{5})},$$
\n and so $$c=5\frac{\sqrt{5(5+2\sqrt{5})}}{2(5-\sqrt{5})}=5\frac{\sqrt{(5+2\sqrt{5})}}{2(\sqrt{5}-1)}=\frac{5}{8}(\sqrt{5}+1)\sqrt{5+2\sqrt{5}} =\frac{5}{8}\sqrt{50+22\sqrt{5}}.$$

For $h\in [h_s,h'_s]$, $A_{min}(h)$ is the area of a regular pentagon minus the area of five equilateral triangles. Each of these regions have sides that depend linearly on $h$. Hence, the formula for  $A_{min}(h)$ must be 
$$A_{min}(h)=\sqrt{5(5+2\sqrt{5})}+c'(h-h_s)(h'_s-h)$$
\n for some constant $c'$ that can be determined by the value $A_{min}(h)$ takes in the middle. If $h=H_{min}/2$, the section is a regular decagon. It is known that the area of such a regular decagon with sides of length $s$ is ${\cal A}=\frac{5}{2} 
\sqrt{5+2\sqrt{5}}s^2.$ Then the equation for $c'$ is
$$\sqrt{5(5+2\sqrt{5})}+c'\frac{1}{20}\cdot 2(5+\sqrt{5})=\frac{5}{2}\sqrt{5+2\sqrt{5}}.$$
So, we get $$c'=5\frac{(\sqrt{5}-2)\sqrt{5+2\sqrt{5}}}{\sqrt{5}+1}=\frac{5}{4}(7-3\sqrt{5})\sqrt{5+2\sqrt{5}}. $$

\n Using this and Mathematica, we obtain 

$$\begin{array}{lr} T_{min}=\frac{1}{15K}\left[ 2 \sqrt{5} \sqrt[4]{2(33112325-14587199 \sqrt{5})}+3 \sqrt[4]{2} (5+\sqrt{5})^{5/4}+\right.
\\ \\ 
\left. 2 \sqrt{10\sqrt{34403829358 \sqrt{5}+76929359725}-1960000-877600 \sqrt{5}}\right].\end{array}
$$

To get the maximum time, Figure~7 (c), we have $H_{max}=2r=\frac{\sqrt{5}+3}{\sqrt{3}}.$ We have to find the spliting points as before:  $h_s$ and $h'_s$ ($h_s<h'_s$). These can be determined from the dihedral angle, $\alpha=\arccos(-\frac{\sqrt{5}}{3})$. We have $\sin(\alpha)=\frac{h_s}{\sqrt{3}}$ which implies $h_s=\frac{2}{\sqrt{3}}$. 
The ratio $\frac{h_s}{H_{max}}=\frac{3-\sqrt{5}}{2}=\frac{1}{2}-\frac{\sqrt{5}-2}{2}$ implies that  $\frac{h'_s}{H_{max}}=\frac{1}{2}+\frac{\sqrt{5}-2}{2}$.

$$\underset{Figure\ 8(a),\  First \ intersection}{\epsfig{file=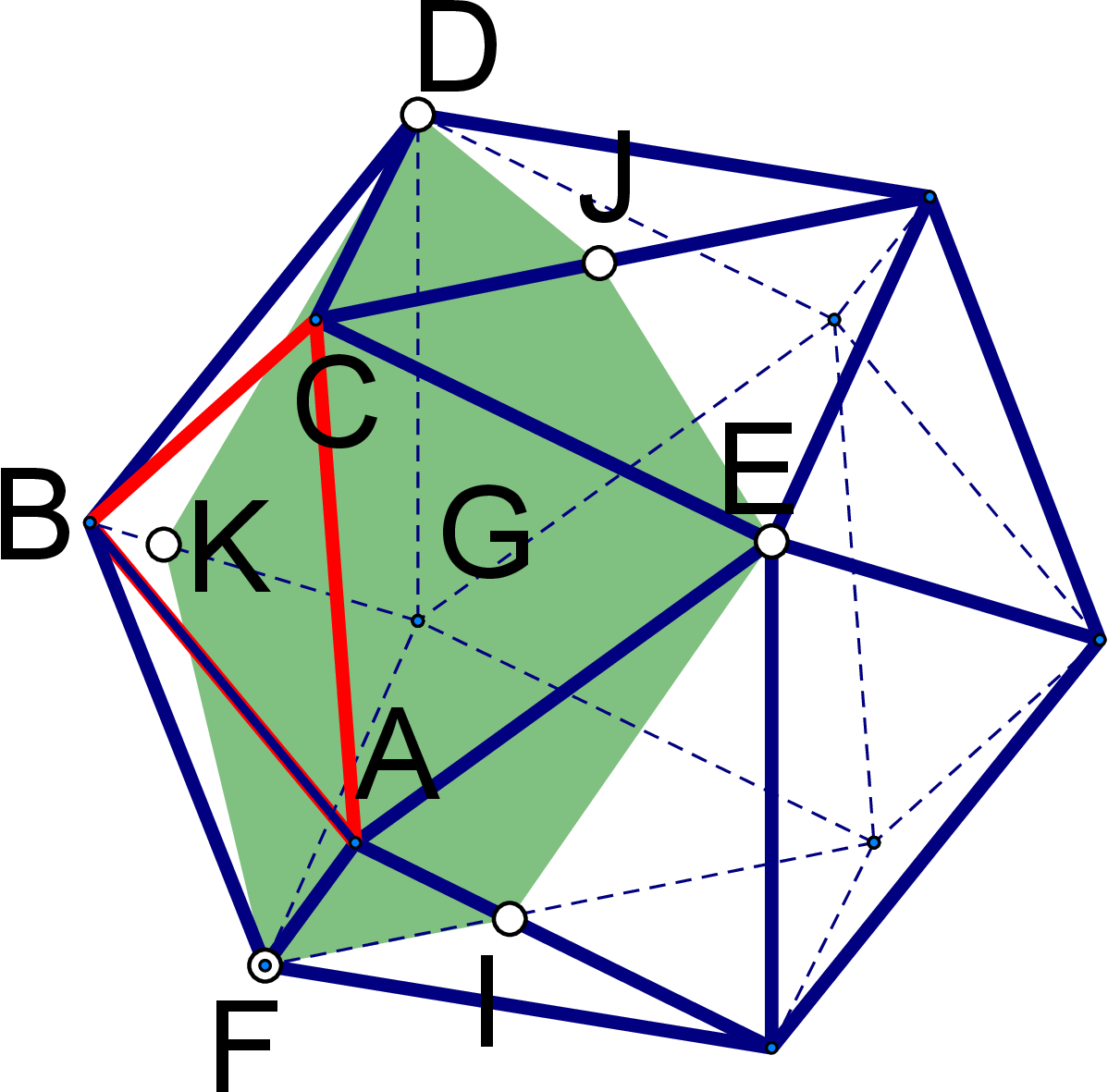,height=2in,width=2in}}$$

We have already observed that the plane of equation $x\phi+z(\phi-1)=\phi^2$ is the plane determined by the points $A$, $B$ adn $C$. We will make the calculations based on Figure~8(a), which is Paciolli's construction which has the advantage of telling us the coordinates of the vertieces. To make the connection with Figure~7(c), we think that the icosahedron in Figure~8 is rotated so that the triangle $ABC$ is horizontal.

 \section{The turn-up number for a solid and a direction} 

To find  $T_{min}$ one may expect to have the solid in a position that has the biggest potential energy of the liquid:

$$P_w=g\delta \int_0^{H_w}A(y)y dy,$$
\n where $w$ is a vector that gives the vertical direction for that solid's position and $\delta$ is the density of the liquid. We will see that this is not the case. This potential energy turns into kinetic energy at the outlet. This is basically Bernoulli’s principle.  For a tiny interval in time, $dt$, this kinetic energy  is $\frac{a(vdt\delta) v^2}{2}$ and using Torricelly's law, this can be written as

$$\delta\frac{a2g\sqrt{2g}h^{3/2}}{2}dt=\delta ga\sqrt{2g}h^{3/2}dt.$$

Integrating over $[0,T]$ and simplifying by $g$ and $\delta$, gives

\begin{equation}\label{eqPotential}
\int_0^{H_w}A(y)y dy=C \int_0^T h_w(t)^{3/2}dt ,
\end{equation}
\n for the constant $C=a\sqrt{2g}=K/k$. If the functions $h_w$, on the common domain, satisfied some inequality like $h_{w_1}(t)>h_{w_2}(t)$  for the same potential energy then, this identity would justify that ``intuition" but in general this is hardly the case. If we multiply (\ref{eq1}) by $h$ and integrate with respect to time over the interval $[0,T]$ we obtain (\ref{eqPotential}) where $C$ is replaced by $K$ as defined earlier. This argument shows that  (\ref{eq1}) is equivalent to Torricelli's Law modulo the constant $k$ which we cannot disregard in reality but we may assume the same for the same solid and get rid of it when comparing times.

These positions, that give the  $T_{min}$ or $T_{max}$,  although well defined from a mathematical point of view, but not necessarily unique, may be difficult to discover. For example, even if  $S$ is a box, say $S=[0, a]\times [0,b]\times [0,c]$ with $a\le b\le c$ gives $T_{min}$ or one may want to take the longest diagonal and push the box so that becomes vertical like in Figure~3 below? 

$$\underset{Figure\ 3}{\epsfig{file=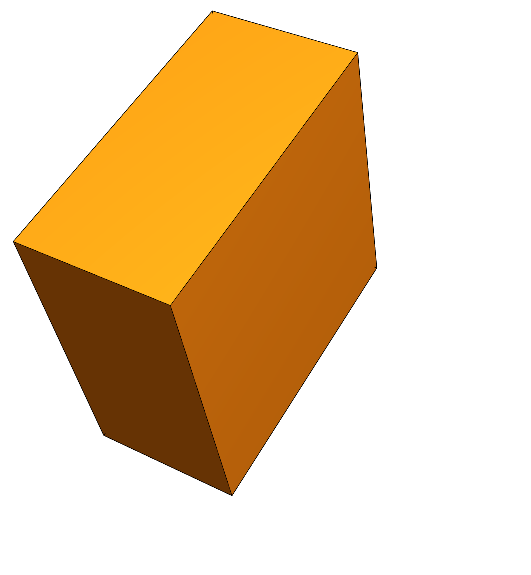,height=3in,width=2in}}$$

So, our calculations are going to be based on a direction $\ell$ ({\it vertical segment} or say {\it vertical line})  picked in advance for a solid in such a way the two planes perpendicular to $\ell$ that are tangent to the solid and contain the solid in between the two are touching the solid at regions that are convex. The last condition is necessary to ensure the drainage of all the liquid. Let us assume that this is happening for every cross-section in between the two planes. This line connects the hole at the lowest point to another point on the boundary of $S$ at the other end of the vertical segment. We can turn the solid upside down and keep the vertical line. 

The ratio $\rho$ from above (\ref{eqMAIN})  is going to be restricted to only the two of these positions relative to the vertical line chosen, and we will use the notation $\rho_{\ell}$ for this and call it the 
{\bf turn-up number}  for a solid $S$ and vertical line $\ell$. For example, if the solid $S$ is a cone, and line $\ell$ is the axis of symmetry,  then the two cases look like the cones in Figure~3 below:

$$\underset{Figure\ 3}{\epsfig{file=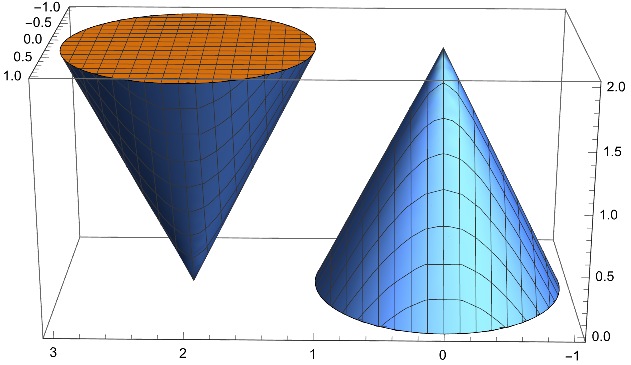,height=2in,width=3in}}$$

  We will show that the cone on the right has a faster drainage time.  In general, we can observe that $H_{max}=H_{min}$ and this value is the distance between the two planes mentioned above. 

Before we do some calculations, let us make another assumption to simplify things. Our solid $S$ is a solid of revolution around $y$-axis of the region between $y$-axis, the graph of $y=f(x)$, $x\in [0,1]$ and $y=1$ having $f$ a strictly increasing function, continuous with $f(0)=0$ and $f(1)=1$. In the cone case, we have $f(x)=x$, for $x\in [0,1]$ for the one on the right in Figure~3. One of the integral involved in  (\ref{eqMAIN}) is then 
$\pi \int_0^1(s^2)^2dy=\frac{\pi}{5}$. The formula for $A(h)$ is $\pi (1-h)^2$ in the opposite position and so
$$\pi \int_0^1(1-s^2)^2ds=\pi(1-2\cdot\frac{1}{3}+\frac{1}{5})=\frac{8\pi}{15}$$
\n which gives $\rho_{\ell}=\frac{8}{15}\cdot (\frac{1}{5})^{-1}= \frac{8}{3}$ for the cone with respect to its axis of symmetry.

Let us observe that in the first case, in general, $A(h)=\pi x^2=\pi (f^{-1}(h))^2$, and for the
upsidedown case, $y=f(x)$ is replaced by $y=1-f(x)$  (we do a reflection across the $x$-axis and a translation upwards by $1$). This means $A(h)=\pi x^2=\pi(f^{-1}(1-h))^2$. Therefore, we have

\begin{equation}\label{eqMAIN2}\rho_{\ell}= \frac{T_{max}}{T_{min}}= \frac{\int_0^1 g(1-s^2)ds}{\int_0^1 g(s^2)ds}
\end{equation} 
\n where $g=(f^{-1})^2$. We just showed that  $\rho_{\ell}= \frac{8}{3}$ if $g(h)=h^2$ for  
 $h\in [0,1]$. Next, we prove that $\rho_{\ell}>1$ for a general class of functions $g$. 

\begin{theorem}
 Let us assume that $g$ is a continuous, non-decreasing function on $[0,1]$ with the property $g(1)>g(0)$. The following inequality takes place
\begin{equation}\label{given} 
\int_0^1g(1-t^2)dt>\int_0^1 g(t^2)dt
\end{equation}
\end{theorem}

\n {\bf Proof:}  We change (\ref{given}) equivalently to 

\begin{equation}\label{given2} 
I:=\int_0^1[g(1-t^2)-g(t^2)]dt>0.
\end{equation}
\n In  (\ref{given2}) we split the integral into two pieces, one over the interval $[0,a]$ and the second over $[a,1]$, where $a=\frac{1}{\sqrt{2}}$:

$$I=\underset{I_1}{\underbrace{\int_0^a[g(1-t^2)-g(t^2)]dt}}+\underset{I_2}{\underbrace{\int_a^1[g(1-t^2)-g(t^2)]dt}}=I_1+I_2.$$

\n In the second integral, we make the substitution $1-t^2=s^2$, $s\ge 0$, 
$$I_2=\int_a^0 [g(s^2)-g(1-s^2)](-\frac{sds}{\sqrt{1-s^2}})=-\int_0^a [g(1-s^2)-g(s^2)]\frac{sds}{\sqrt{1-s^2}}.$$
\n With this calculation, we go back to  (\ref{given2}) and continue:
  
$$I=I_1+I_2=\int_0^a[g(1-t^2)-g(t^2)]dt-\int_0^a [g(1-s^2)-g(s^2)]\frac{sds}{\sqrt{1-s^2}}\implies 
$$

$$ I=\int_0^a[g(1-t^2)-g(t^2)]\left(1-\frac{t}{\sqrt{1-t^2}}\right) dt .$$

\n Clearly, to prove  (\ref{given2}) it suffices to show that the integrant is non-negative and strictly greater than zero for some subinterval of $[0,a]$. For every $t\in[0,a]$ we have $1-t^2\ge t^2$, and since $h$ is non-decreasing  $g(1-t^2)\ge g(t^2)$. This inequality is strict for $t$ small by the hypothesis of our problem and the continuity of $g$. Also, $1-\frac{t}{\sqrt{1-t^2}}\ge 0$ is equivalent to  $\sqrt{1-t^2}\ge t$
which we already mentioned in the form $1-t^2\ge t^2$. Hence,  the integral $I$ is a positive real number.\eproof

\par \vspace{0.1in}

 Taking $g(x)=x^{1/n}$ in  (\ref{given})  we obtain the following inequality. 
\begin{corollary}   For every $n\in \mathbb N$, we have 
\begin{equation}\label{given3}
\int_0^1(1-t^2)^{\frac{1}{n}}dt > \frac{n}{n+2}.
\end{equation}
\end{corollary}

 \section{Turn-up number equal to 1} 

In this section, we are still considering a solid of revolution around the y-axis but the region is determined by the y-axis, and the graph of a continuous function implicitly given by $x^2=g(y)$, $y\in [0,1]$ and $g(0)=g(1)=0$. Certainly, we want to obtain a solid $S$ that has no singular points, so let us assume also that $g(y)>0$ for all $y\in (0,1)$ (as in Figure~4).

$$\underset{Figure\ 4,\  x=\sqrt{g(y)}}{\epsfig{file=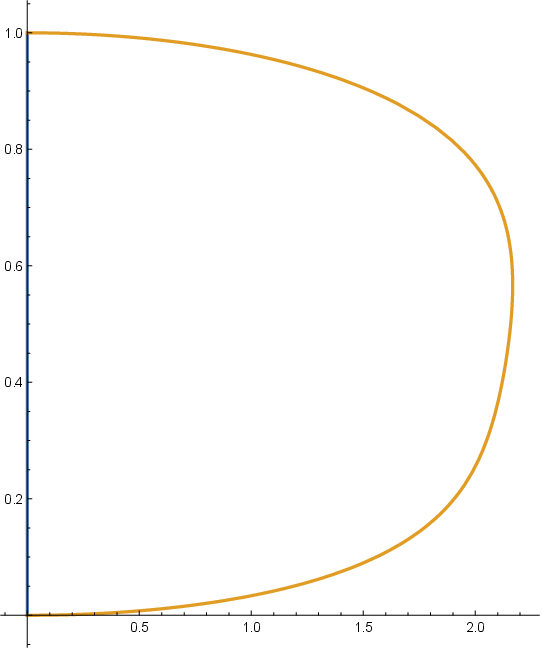,height=2.5in,width=2.5in}}$$ 

Clearly, if $g$ is symmetric with respect to the horizontal line $y=1/2$, i.e., $g(y)=g(1-y)$ for all $y$, we obtain $\rho=1$, but we are interested in solids that are not symmetric. The first step is to try a polynomial in $y$ and $1-y$. For example, $y(1-y)$ satisfies the requirements but it is symmetric, and the same is true for  $y^m(1-y)^m$, which $m\ge 1$. So, we may want to try a combination of the form $g(y)=Cy^m(1-y)^n+y^p(1-y)^q$ for some constants $C$, $m$, $n$, $p$ and $q$ and determine $C$ so that $\int_0^1 g(s^2)ds=\int_0^1 g(1-s^2)ds$. It turns out that one of the smallest such examples is given by 
$$g(y)=29y^2(1 - y) + 33(1 - y)^4y$$
and one can check that $\int_0^1 g(s^2)ds=\frac{302}{105}=\int_0^1 g(1-s^2)ds$. The graph of $x=\sqrt{g(y)}$ is actually in Figure~4. This implies $\rho=1$, and we obtain the solid in  Figure~5 below:

$$\underset{Figure\ 5(a),\  g(1-z)=x^2+y^2}{\epsfig{file=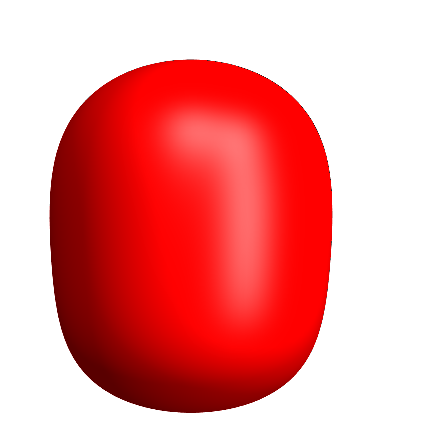,height=4in,width=3in}}\underset{Figure\ 5(b),\  g(z)=x^2+y^2 }{\epsfig{file=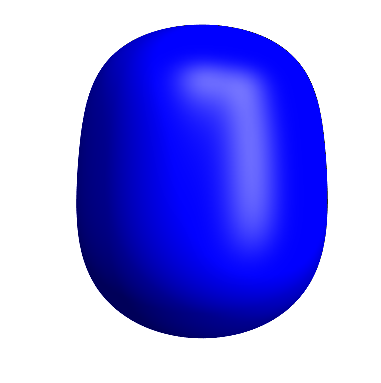,height=4in,width=3in}}$$
The function $g$ satisfies the requirements we wanted and in addition, 
it can be shown that $g$ is concave downwards. Also, one can check that $g'(0)=33$ and $g'(1)=-29$. Then the solid of revolution obtained is convex, having a tangent plane at any point on its boundary. This is because we are rotating $x=\sqrt{g(y)}$ and so $\frac{dy}{dx}=\frac{\sqrt{g(y)}}{g'(y)}$ which at $y=0$ and $y=1$ is zero. Now, if we compute the potential energy (divided by $\pi g\delta$, $\delta$ the density of the liquid) for the two cases we get 

$$ \int_0^1 g(y)ydy=\frac{247}{140}, \ \text{and}\  \   \int_0^1 g(1-y)ydy=\frac{184}{105}, \ \text{with} \ \frac{247}{140}-\frac{184}{105}=\frac{1}{84}. $$

This is to say that the red solid (in Figure~5) has more potential energy than the blue one, as expected.  
One can find infinitely many polynomial functions like this one.  Here is another example constructed with the same idea: 
$$G(y)= 13y^2(1 - y)^6 + 9y^3(1 - y)^2,\ y\in [0,1]$$
\n which has the property that at the points on the $y$-axis, there are pointing ends as in Figure~6 (not smooth):
$$\underset{Figure\ 6,\   g(z)=x^2+y^2 \ \ g(1-z)=(x-2)^2+(y-2)^2}{\epsfig{file=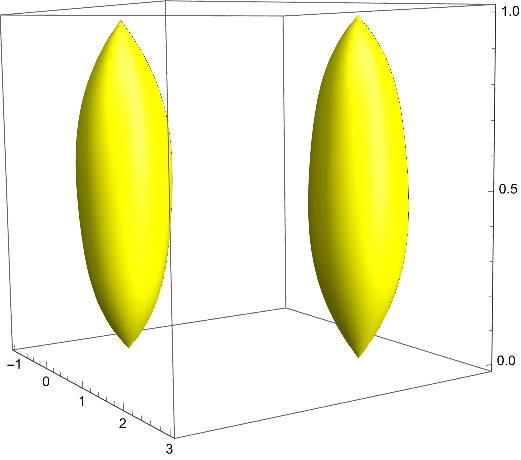,height=2in,width=3in}}$$

\subsection{Characterization of $\rho=1$} 
We assume that $x=g(y)$ is continuous and defined for all $y\in [0,1]$. The turn-up number is associated with the solid obtained by rotating
the graph of $g$ around the y-axis and bounded from above by the plane $y=1$ and from below by $y=0$. For computing $\rho$ we use the axis of rotation which is the y-axis. We have seen that $\rho=1$ is equivalent to $\int_0^1 g(s^2)ds=\int_0^1 g(1-s^2)ds$ or $\int_0^1 [g(1-s^2)-g(s^2)]ds=0$.

Let us make the substitution $s=\sin t$: 
$$0=\int_0^{\pi/2} [g(\cos ^2t)-g(\sin ^2 t)]\cos t dt=\int_0^{\pi/4}... + \int_{\pi/4}^{\pi/2}...$$
In the last integral, we make another substitution, namely $u=\pi/2-t$, and obtain 
$$\int_{\pi/4}^{\pi/2}  [g(\cos ^2t)-g(\sin ^2 t)]\cos t dt=\int_{\pi/4}^{0}  [g(\sin ^2t)-g(\cos ^2 t)]\sin t dt.$$
We see that then $\rho=1$ is equivalent to

$$\int_0^{\pi/4} [g(\cos ^2t)-g(\sin ^2 t)](\cos t-\sin t) dt=0.$$

\n We observe that $\cos t-\sin t=\sqrt{2}\cos(\frac{\pi}{4}+t)$. Hence we have the following characterization.

\begin{theorem}
 Let us assume that $g$ is as described above. Then $\rho=1$ is equivalent to
\begin{equation}\label{eqeq} 
 \int_0^{\pi/4} [g(\cos ^2t)-g(\sin ^2 t)]\cos(\frac{\pi}{4}+t)dt=0.
\end{equation}
\end{theorem}
So, in order to construct some $g$ that is not symmetric all we have to do is define $g$ on $[0,1/2]$ arbitrarily (just $g(y)>0$ for $y\in (0,1/2]$ ) and then extend it 
by $g(\cos ^2t)=g(\sin ^2 t)+h(t)$, $t\in [0,\pi/4]$  for some non-zero continuous function $h$ so that  $\int_0^{\pi/4} h(t)\cos(\frac{\pi}{4}+t)dt=0.$ We would like to have $g$ a continuous function at $y=1/2$ so $h(\pi/4)=0$. Also, we need to have $g(y)>0$ for all $y\in (1/2,1)$. This is equivalent to 
$h(t)>-g(\sin ^2 t)$ for all $t\in (0,\pi/4)$.
\par\vspace{0.2in}

For more developments on this topic we refer to the thesis in \cite{bianca}. 

\end{document}